\documentclass[journal,10pt, letterpaper]{IEEEtran}  

\IEEEoverridecommandlockouts                              

\usepackage{xcolor}
\definecolor{update}{RGB}{0,160,255}
\usepackage{algorithm}
\usepackage{algorithmic}

\def\F{\mathcal{F}}
\def\U{\mathcal{U}}
\def\X{\mathcal{X}}

\def\powip{^{(i+1)}}
\def\powi{^{(i)}}

\usepackage{graphics} 
\usepackage{epsfig} 
\usepackage{times} 
\usepackage{amsmath} 
\usepackage{amssymb}  
\usepackage{amsfonts}

\def\vec{\mathrm{vec} \, }

\usepackage{tikz}
\usepackage{pgfplots} 
\pgfplotsset{compat=1.17}
\usepgfplotslibrary{external}
\usetikzlibrary{positioning}
\usepackage{comment}
\usepackage{subcaption}

\newtheorem{lemma}{Lemma}[section]
\newtheorem{assumption}{Assumption}[section]

\DeclareMathOperator{\argmin}{argmin}
\allowdisplaybreaks[4]

\title{Scenario Model Predictive Control for Data-based Energy Management in Plug-in Hybrid Electric Vehicles
}

\author{Sebastian East$^{1}$ and Mark Cannon$^{2}$
\thanks{This work was supported by the Engineering and Physical Sciences Research Council}
\thanks{$^{1}$Sebastian East is with the Department of Aerospace Engineering, University of Bristol,
        Bristol, BS8 1TR, UK
        {\tt\small sebastian.east@bristol.ac.uk}}%
\thanks{$^{2}$Mark Cannon is with the Department of Engineering Science, University of Oxford,
        Oxford, OX1 3JP, UK
        {\tt\small mark.cannon@eng.ox.ac.uk}}%
}

\begin{document}

\maketitle
\pagestyle{plain}

\begin{abstract}

One of the major limitations of optimization-based strategies for allocating the power flow in hybrid powertrains is that they rely on predictions of future power demand. These predictions are inherently uncertain as they are dependent on complex human behaviours that are challenging to model accurately. This paper proposes a data-based scenario model predictive control framework, where the inputs are determined at each control update by optimizing the power allocation over multiple previous examples of a route being driven. The proposed energy management optimization is convex, and results from scenario optimization are used to bound the confidence that the one-step-ahead optimization will be feasible with given probability. It is shown through numerical simulation that scenario model predictive control (MPC) obtains the same reduction in fuel consumption as nominal MPC with full preview of future driver behaviour, and that the scenario MPC optimization can be solved efficiently using a tailored optimization algorithm.  

\end{abstract}

\section{INTRODUCTION}

Plug-in hybrid electric vehicles (PHEVs) have a powertrain that includes an internal combustion engine, an electric motor, and a battery that can be charged from an external source. This configuration is typically chosen to obtain the benefits of all-electric propulsion (e.g. zero tailpipe emissions and regenerative braking) with a smaller battery than would be required for an all-electric vehicle, whilst also ensuring that the total combined range available before a recharge is at a useful level. The key differentiators between a PHEV and a `mild' hybrid electric vehicle (HEV) are the size of the battery and the charging method: the battery in a HEV is much smaller and cannot be charged from an external source, and is therefore typically only used for regenerative braking and to accelerate the vehicle from stationary.

Multiple studies have demonstrated that the the energy consumption (and consequently, the range) of a PHEV is strongly determined by the controller that allocates the driver's power demand to the engine and motor (a comprehensive review of control methods is provided in \cite{Martinez}). A simple strategy is to allocate all of the power demand to the motor until the battery is empty, and then only use the motor to regeneratively brake and start the vehicle from a standstill until the end of the journey. However, blended strategies, where power is delivered from both sources simultaneously and the fraction delivered from each is varied throughout the journey, are potentially much more effective. In particular, optimization-based strategies, where the power allocation is determined using numerical optimization of a mathematical representation of the powertrain can considerably reduce fuel consumption without affecting total battery discharge. Many approaches to optimization-based energy management have been considered, including dynamic programming \cite{Moura}, Pontryagin's minimum principle \cite{Stockar}, equivalent consumption minimization strategy \cite{Sciarretta}, and model predictive control (MPC) \cite{Sun}. 

A major limitation of any optimization-based strategy is that its effectiveness is conditional on an accurate prediction of the future behaviour of the driver (specifically, the velocity of the vehicle and the torque at the wheels throughout the remaining journey). This prediction is non-trivial, and depends on a range of factors that may be challenging to model and predict, including the route, the geometry of the road, traffic conditions, and driving style. A common approach is to fit a model (such as a neural network \cite{Sun2, Chen2014} or Markov-based model \cite{Zeng, Cairano2014}) to driver behaviour data, and then use the model to generate predictions that are fed into the energy management optimization. One of the fundamental limitations of this approach is the model must necessarily be conditioned on a large number of environmental variables (e.g. road geometry, speed limits) that vary throughout a given journey and may themselves be uncertain. Consequently, these techniques can only generally be used to generate `short-horizon' predictions (e.g.\ up to $\sim$20\,s into the future); long-horizon predictions (i.e. predictions of an entire journey lasting minutes or hours) remain an open problem \cite{Zhou2019}.

In this paper, a data-based scenario MPC is proposed where the power-split is determined at each sampling instant by optimizing the predicted fuel consumption over multiple previous examples of the route being driven. The primary benefit of this approach is that predictions of driver behaviour (that inherently account for route-specific variables, such as road geometry and speed limits) can be made over arbitrarily long horizons, and an explicit driver model is not required. The proposed method is applicable both when a vehicle repeatedly completes the same route (e.g. public transportation and delivery/collection vehicles) and when the route to be driven is provided by a navigation device (or the vehicle is being driven autonomously). It could also be combined with a statistical method that automatically identifies the route being driven from previous driver behaviour (such as that proposed in \cite{Larsson}), but this is not investigated in this paper.

The main contributions of this paper are as follows:
\begin{enumerate}
\item A convex scenario MPC framework that minimizes fuel consumption across multiple predictions of future driver behaviour. The approach incorporates nonlinear losses in the engine, motor, and battery; limits on the energy stored in the battery; and limits on the power delivered by both the electric motor and the internal combustion engine. Importantly, the scenario MPC framework directly optimizes the state-of-charge trajectory over the entire prediction horizon, and does \emph{not} rely on tracking a reference trajectory (as proposed in, for example, \cite{Sun}). Results from scenario optimization \cite{Schildbach2013} provide a bound on the confidence that the one-step-ahead MPC optimization will be feasible with a given probability, given the number of scenarios considered in the optimization. 
\item An alternating direction method of multipliers (ADMM) algorithm  tailored to the structure of the scenario MPC optimization problem, where each iteration has $\mathcal{O}(NS)$ memory and computational cost (where $N$ is the prediction horizon length and $S$ is the number of scenarios).
\item A set of numerical studies in which a nominal MPC strategy with access to accurate predictions of future driver behaviour is compared against scenario MPC, for which predictions are generated directly from previously recorded journeys on the same route (i.e. the controller has no preview of the journey being simulated). Remarkably, it is found that the scenario MPC obtains the same fuel consumption as nominal MPC. 
\end{enumerate}

This is not the first paper proposing scenario MPC for energy management: it was previously applied to a charge sustaining hybrid-electric vehicle in \cite{Josevski}, but that paper considered a linearized vehicle model, a prediction horizon of only 10 samples, and scenarios were generated using perturbations of a nominal prediction obtained from the non-parameteric model proposed in \cite{Josevski2}.

This paper builds on the authors' previous work on convex optimization-based energy management in \cite{CDC2018} and \cite{EastTCST}. In those cases only a single prediction of future driver behaviour was considered, and the extension to multiple scenarios presented here is original work. Additionally, the authors have previously investigated the use of graphics processing units (GPUs) for accelerating ADMM when used to solve similar scenario-based resource allocation problems in \cite{Qureshi2019}. PHEV energy management was also used as an illustrative example in that paper, but the scenarios were generated using random perturbations of a nominal trajectory (rather than the data-based framework proposed in this paper), and the energy management problem was considerably simplified, and therefore not suitable for supervisory control as proposed here.

The rest of this paper is structured as follows. In Section \ref{section::problem_formulation} the energy management problem is formulated mathematically. Section \ref{section::scenario_mpc} presents the SMPC controller, and the algorithm used to solve the corresponding optimization problem is presented in Section \ref{section::optimization}. The results of numerical experiments are presented and discussed in Section \ref{section::numerical_experiments}.

\textbf{Notation}
The set of integers from $a \in \mathbb{Z}$ to $b \in \mathbb{Z}$ (inclusive) is denoted by $\{a, \dots, b\}$. For a given set $\mathbb{S}$, the set of non-negative (positive) elements of $\mathbb{S}$ is denoted by $\mathbb{S}_+$ ($\mathbb{S}_{++}$), and the set of non-positive (negative) elements of $\mathbb{S}$ is denoted by $\mathbb{S}_-$ ($\mathbb{S}_{--}$). The Minkowski sum of sets $A,B$ is $A\oplus B =\{x+y: x\in A, \, y\in B\}$.
The notation $M_{:,i}$ ($M_{i,:}$) refers to the $i$th column (row) of matrix $M$. The function $\vec (M)$ returns a column vector containing the columns of matrix $M$ stacked sequentially from left to right. All inequalities are considered element-wise in the context of vectors. The symbols $\mathbf{1}$ and $\mathbf{0}$ respectively denote a vector of ones and a vector of zeros (dimensions can be inferred from their context).

\section{Problem Formulation} \label{section::problem_formulation}

\tikzstyle{input} = [coordinate]
\tikzstyle{sum} = [draw, fill=black!10, circle, node distance=0.5cm]
\tikzstyle{block} = [draw, fill=black!10, rectangle, 
    minimum height=3em, minimum width=3em] 

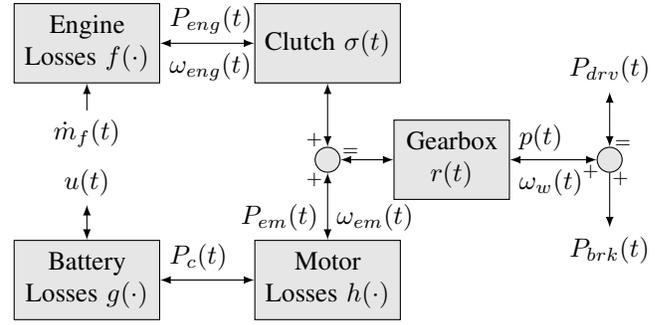
\begin{figure}
\begin{center}
\begin{tikzpicture}[auto, node distance=2.8cm,>=latex]
\node [name = input] {$\dot{m}_f(t)$};
\node [block, name=engine, above of=input, node distance = 1.2cm, text width=1.7cm, align=center] {Engine Losses $f(\cdot)$};
\node [block, name=clutch, right of=engine, node distance = 3.2cm, text width=1.7cm, align=center] {Clutch $\sigma(t)$};
\node [sum, below of=clutch, node distance = 1.55cm, name=sum] {};
\node [block, right of = sum, name=gearbox, node distance=1.65cm, text width=1.3cm, align=center] {Gearbox $r(t)$};
\node [block, below of = sum, node distance = 1.60cm, name=motor, text width=1.7cm, align=center] {Motor Losses $h(\cdot)$};
\node [block, left of=motor, node distance=3.2cm, name=batteryloss, text width=1.7cm, align=center] {Battery Losses $g(\cdot)$};
\node[above of=batteryloss, node distance=1.3cm, name=supercapinput] {$u(t)$};
\node [sum, right of=gearbox, node distance = 2.1cm, name=sum2] {};
\node [below of =sum2, node distance=1.2cm, name=brake] {$P_{brk}(t)$};
\node [above of =sum2, node distance=1.2cm, name=drive] {$P_{drv}(t)$};
\draw [->] (input) -- (engine);
\draw [<->] (engine) -- (clutch) node [pos=0, below right] {$\omega_{eng}(t)$} node [pos=0, above right] {$P_{eng}(t)$};
\draw [<->] (clutch) -- (sum);
\draw [<->] (motor) -- (sum) node [pos=0, above right] {$\omega_{em}(t)$} node [pos=0, above left] {$P_{em}(t)$};
\draw [<->] (gearbox) -- (sum) node [midway, above] {};
\draw [<->] (batteryloss) -- (motor) node [pos=0, above right] {$P_{c}(t)$};
\draw [<->] (supercapinput) -- (batteryloss);
\draw [<->] (gearbox) -- (sum2) node [pos = 0, above right] {$p(t)$} node [pos = 0, below right] {$\omega_w(t)$};
\draw [->] (sum2) -- (brake);
\draw [<->] (sum2) -- (drive);
\node [below right = -0.09cm and -0.21cm of sum2] { \scriptsize{$+$}};
\node [below left = -0.18cm and -0.10cm of sum2] { \scriptsize{$+$}};
\node [above right = -0.075cm and -0.18cm of sum2] { \scriptsize{$=$}};
\node [above left = -0.075cm and -0.18cm of sum] { \scriptsize{$+$}};
\node [below left = -0.075cm and -0.18cm of sum] { \scriptsize{$+$}};
\node [above right = -0.18cm and -0.05cm of sum] { \scriptsize{$=$}};
\end{tikzpicture}
\bigskip
\caption{Diagram of the loss functions and variables that are used in the mathematical formulation of the energy management problem. The directions of the arrows represents the direction of power transfer. Symbols are included to indicate how the powerflows are combined at the summing junctions.}
\label{figure::PHEV_powerflow_diagram}
\end{center}
\end{figure}

This section presents a mathematical model of a parallel, pre-transmission PHEV powertrain that is used to formalize the energy management problem, and later used for simulation. The model used in this paper is the same as that used in \cite{EastTCST} (and is also similar to vehicle models used in the PHEV energy management literature, e.g. \cite{Moura, Stockar,Sciarretta,Sun}), and the reader is therefore referred to \cite{EastTCST} for a more in-depth exposition; only the details required to understand the following sections are included here. Figure \ref{figure::PHEV_powerflow_diagram} shows a simplified diagram of the powertrain under consideration, for which the relevant behaviour of the driver can be described completely by the vehicle's velocity, $v(t)$, and the road gradient, $\theta(t)$, for all $t \in [0,T]$, where $T \in (0, \infty)$ is the duration of the journey (the problem is \emph{not} to control the behaviour of the vehicle in an autonomous sense). 
From this description of driver behaviour, the rotational speed of the wheels, $\omega_w(t)$, and power demanded at the wheels, $P_{drv}(t)$, are defined by 
\begin{equation}\label{equation::wheel_speed}
\omega_{w}(t) := v(t) / r_w 
\end{equation}
(where $r_w \in \mathbb{R}_{++}$ is the effective radius of the wheels) and the longitudinal vehicle model
\newpage
\begin{multline}\label{equation::longitudinal_model}
    P_{drv}(t) := \Bigg[m\dot{v}(t) + \frac{1}{2}\rho_a v(t)^2C_dA \\
    + C_rmg \cos \theta(t) + mg \sin \theta(t) \Bigg] v(t),
\end{multline}
where $m \in \mathbb{R}_{++}$ is the mass, $\rho_a \in \mathbb{R}_{++}$ is the density of air, $C_d \in \mathbb{R}_{++}$ is the drag coefficient, $A \in \mathbb{R}_{++}$ is the frontal area, $C_r\in \mathbb{R}_{++}$ is the rolling resistance, and $g \in \mathbb{R}_{++}$ is the acceleration due to gravity (it is also assumed that $v(t)$ is differentiable everywhere). Consequently, there are four parameters that must controlled for all $t \in [0,T]$ in order that the operation of the powertrain is well-defined: 
\begin{enumerate}
\item The state of the engine, $\sigma(t) \in \{0, 1\}$, where $\sigma(t) = 1$ implies that the engine is on and the clutch engaged, and $\sigma(t) = 0$ implies that the engine is off and the clutch disengaged.
\item The gear selection, $r(t) \in \mathcal{R}$, where $\mathcal{R} \subset \mathbb{R}$ is the set of available gear ratios.
\item The mechanical braking power, $P_{brk}(t) \in \mathbb{R}_{-}$.
\item The power delivered by the engine, $P_{eng}(t) \in \mathbb{R}_{+}$ (engine braking is not considered). 
\end{enumerate}

The power delivered from the motor, $P_{em}(t) \in \mathbb{R}$ does not need to be explicitly selected as it is determined implicitly by the other four specified control variables (the value of $P_{brk}(t)$ determines the relative fraction of braking power delivered by the motor and mechanical brakes).

Once the gear selection and engine state have been determined, the rotational speeds of the engine, $\omega_{eng}(t) \in \mathbb{R}_{+}$, and motor, $\omega_{em}(t)$, are defined by
\begin{equation}\label{equation::rotational_speed}
\left[ \omega_{eng}(t), \omega_{em}(t) \right] := \begin{cases} \left[ r(t)\omega_{w}(t), r(t) \omega_{w}(t) \right] & \sigma(t) = 1 \\
\left[ 0, r(t) \omega_w(t) \right] & \sigma (t) = 0 \end{cases}
\end{equation}
for all $t  \in [0, T]$. Once the braking power has been determined, the power, $p(t)$, delivered by the drivetrain to the axles is 
\begin{equation}\label{eqn::pwtn_output}
p(t) := P_{drv}(t) - P_{brk}(t) \quad \forall t \in [0,T].
\end{equation}
The drivetrain components are assumed to be 100\% efficient, so once the engine output power, $P_{eng}(t)$, has been determined, the power output required of the motor is defined by
$$
P_{em}(t) := \begin{cases} p(t) - P_{eng}(t) & \sigma(t) = 1 \\ p(t) & \sigma(t) = 0 \end{cases}
$$
for all $t \in [0, T]$.

The engine's rate of fuel consumption, $\dot{m}_f(t)$, is modelled as a static input-output map, $f: \mathbb{R} \times \mathbb{R} \mapsto \mathbb{R}_+$, defined by
$$
\dot{m}_f(t) := f( P_{eng}(t), \omega_{eng} (t) ).
$$
Similarly, the motor's output power, $P_c(t)$, is modelled as a static input-output map, $h: \mathbb{R} \times \mathbb{R} \mapsto \mathbb{R}$, defined by
$$P_{c} (t) := h(P_{em}(t), \omega_{em}(t) ).
$$
The battery is modelled as an equivalent circuit consisting of an internal resistance only, so that the rate of depletion of the battery's internal energy, $u(t)$, is modelled by the function $g: \mathbb{R} \times \mathbb{R}_{++} \times \mathbb{R}_{++} \mapsto \mathbb{R}$, defined by
\begin{align*}
u (t) & := g(P_c(t), V(t), R(t)) \\
& : =  \frac{V(t)^2}{2 R(t)}\left(1 -   \sqrt{1 - \frac{4 R(t)}{V(t)} P_c(t)} \right)
\end{align*}
$\forall P_c(t) \leq \frac{V(t)}{4R(t)}$, where $V(t) \in \mathbb{R}_{++}$ and $R(t) \in \mathbb{R}_{++}$ are the battery's open circuit voltage and internal resistance. Consequently, the internal energy (or, synonymously, the state-of-charge) of the battery is defined by
$$
E(t) = E(0) - \int_0^t u(t) \, \mathrm{d} t \quad \forall t \in [0, T]. 
$$
For the purposes of this formulation it is assumed that $u(t)$ is Riemann integrable. The objective of optimization-based energy management is to minimize the fuel consumed during a given journey, i.e. to obtain values $d(t)^\star := (\sigma(t)^\star, r(t)^\star, P_{brk}(t)^\star, P_{eng}(t)^\star)$ that minimize
\begin{equation}\label{equation::optimal_contrinuous}
\inf_{
\substack{
d(t) \in \mathcal{D}(t) \\
\forall t \in [0, T]
} 
} m_f(T).
\end{equation}
Here $\mathcal{D}(t)$ is a time-dependent set that enforces hard constraints on the torque, rotational speed, power, and energy of each of the components illustrated in Figure \ref{figure::PHEV_powerflow_diagram} (see \cite{EastTCST} for a more complete description). In particular, the engine power has instantaneous upper and lower bounds (which also implicitly impose bounds on the motor and battery power), 
$$
\underline{P}_{eng}(t) \leq P_{eng}(t) \leq \overline{P}_{eng}(t) \quad \forall t \in [0, T]
$$
and the state-of-charge of the battery is constrained between upper and lower limits, 
$$
\underline{E} \leq E(t) \leq \overline{E} \quad \forall t \in [0, T ].
$$
Obtaining a solution to (\ref{equation::optimal_contrinuous}) is intractable in general due to the complex constraints on system operation, and the fact that the solution depends on the power and velocity required of the vehicle in the future, which cannot be known with complete precision \emph{a priori}.

\section{Scenario Model Predictive Control}\label{section::scenario_mpc}

We propose a MPC framework for solving
(\ref{equation::optimal_contrinuous}) approximately,
in which the control variables are updated at a frequency of 1\,Hz\footnote{The approach can be readily extended to an arbitrary sampling interval, subject to the associated MPC optimization solution times.} throughout the journey by solving a discrete-time optimal control problem, and held constant between each update. To aid the presentation of the proposed framework, we begin by considering the deterministic case where at each variable update instant the controller has a single prediction of the vehicle's future velocity, $v:=(v_0, \dots, v_{N-1}) \in\mathbb{R}^N$, and the gradient of the road, $\theta:=(\theta_0, \dots, \theta_{N-1})\in\mathbb{R}^N$, where $N \in \mathbb{N}$ is the length of the prediction horizon\footnote{Note that, for the example variable $v$, the notation $v(t)$ is used to represent the simulated velocity at time $t$, whereas the notation $v_k$ is used to represent the velocity \emph{predicted} by the MPC at time $v(t + k)$.}; later this formulation will be extended to consider the case of multiple predictions. It is assumed that the predictions are also sampled at 1 Hz, although this is not required and the methods presented here can be extended to an arbitrary sampling interval. The predicted rotational speed of the wheels, $\omega_{w,k} \in \mathbb{R}$, and driver's power demand, $P_{drv,k} \in \mathbb{R}$, can then be determined for all $k \in \{0, \dots, N-1 \}$ using similar forms of \eqref{equation::wheel_speed} and \eqref{equation::longitudinal_model} (the accelerations $\dot{v}_k$ can be calculated using numerical differentiation of $v$).

\subsection{Convex Model Predictive Control}\label{section::convex_MPC}

It has previously been demonstrated in \cite{CDC2018} that if the gear selection, braking power, and engine switching control are pre-determined (by either a heuristic or optimization-based method), then the problem of determining the optimal power split between the motor and internal combustion engine w.r.t.~\eqref{equation::optimal_contrinuous} can be reduced to a convex optimization problem under minor further approximation\footnote{There is an alternative method for achieving a convex formulation, as presented in \cite{Egardt2014}, but the method presented in \cite{CDC2018} results in an optimization problem with preferable structure. The solutions obtained using both methods are mathematically equivalent.}, and this is the approach taken here. Again, only sufficient detail is reproduced here to understand the following development, and interested readers are referred to \cite{CDC2018} for more details.

To achieve the convex energy management formulation, firstly the engine and motor maps are approximated as strictly convex quadratic functions\footnote{This is a common approach, and has been shown to achieve a very close approximation of the true loss maps (see for example \cite[\S 7]{MURGOVSKI2012106}).}:
\begin{multline*}
f( P_{eng}(t), \omega_{eng} (t) ) \approx \\ 
\alpha_2 (\omega_e (t) ) P_{eng}(t)^2 + \alpha_1 (\omega_e (t) ) P_{eng}(t) + \alpha_0 (\omega_e (t) ),
\end{multline*}
where $\alpha_2(x) > 0$ $\forall x \in \textrm{dom}(\alpha_2 )$, and
\begin{multline*}
h( P_{em}(t), \omega_m (t) ) \approx \\
\beta_2 (\omega_m (t) ) P_{em}(t)^2 + \beta_1 (\omega_m (t) ) P_{em}(t) + \beta_0 (\omega_m(t) ),
\end{multline*}
where $\beta_2(x) > 0$ $\forall x \in \textrm{dom}(\beta_2 )$. Then, at each MPC update when in operation, the following steps are taken to obtain the convex formulation as presented in \cite{CDC2018}:
\begin{enumerate}
\item The predicted values of engine state, $\sigma_k$, braking power, $P_{brk,k}$, and gear selection, $r_k$, are determined using the external control strategy for all $k \in \{ 0, \dots , N-1 \}$.
\item The predicted values of powertrain output power, $p_k$, are determined using a similar form of \eqref{eqn::pwtn_output} (i.e. indexed by $k$ instead of $t$) for all $k \in \{0, \dots, N-1 \}$.
\item The set $\mathcal{P}$ is defined as the subset of values of ${0, \dots, N-1}$ where the engine is switched on, i.e.
$$
\mathcal{P}:= \{k : \sigma_k = 1 \}.
$$
\item The predicted values of engine speed, $\omega_{eng,k}$, and motor speed, $\omega_{em,k}$, are determined for all $k \in \{0, \dots, N-1 \}$ using a similar form of \eqref{equation::rotational_speed}.
\item The predicted rate of fuel consumption, $\dot{m}_{f,k} \in \mathbb{R}_+$, and battery output power, $P_{c,k} \in \mathbb{R}$, are defined by the quadratic functions
\newpage
\begin{align*}
\dot{m}_{f,k} & := f_k (P_{eng,k}) \\
&: = \alpha_{2,k} P_{eng,k}^2 + \alpha_{1,k} P_{eng,k} + \alpha_{0,k} \\
P_{c,k} & : = h_k (P_{em,k}) \\
& : = \beta_{2,k} P_{em,k}^2 + \beta_{1,k} P_{em,k} + \beta_{0,k},
\end{align*}
where $\alpha_{i,k} := \alpha_i (\omega_{eng,k})$ and $\beta_{i,k} := \beta_i (\omega_{em,k}),$ for all $i \in \{0,1,2\}$ and $k \in \{0, \dots, N-1 \}$.
\item The open circuit voltage and resistance of the battery, $V_k$ and $R_k$, are assumed to be fixed \emph{a priori} and not considered as part of the optimization (for small variations in battery SOC it can be effective to set $(V_k, R_k) = (V(t), R(t))$ for all $k \in \{0, \dots, N-1 \}$).
\item The predicted rate of battery depletion, $u_k \in \mathbb{R}$, is defined as equal to the function
\begin{align*}
u_k & := g_k(P_{em,k}) \\
& : = \frac{V_k^2}{2 R_k}\left(1 -   \sqrt{1 - \frac{4 R_k}{V_k} h_k(P_{em,k})} \right)
\end{align*}
for all $k \in \{0, \dots, N-1 \}.$
\item Values of $\underline{u}_k \in \mathbb{R}$ and $\overline{u}_k \in \mathbb{R}$ are generated for all $k \in \{0, \dots, N-1 \}$ so that the box constraints ${\underline{u}_k \leq u_k \leq \overline{u}_k}$ enforce instantaneous power limits on all powertrain subsystems, constrain the arguments of $f_k(\cdot )$ and $h_k(\cdot )$ to the domain in which they are strictly increasing, and deliver sufficient electrical power when the engine is off (i.e. for all $k \notin \mathcal{P}$).
\end{enumerate}

Consequently, the fuel consumption at each timestep $k$ can be modelled in terms of the battery's internal power by the convex function $F_k: [\underline{u}_k, \overline{u}_k] \mapsto \mathbb{R}_{+}$, defined by
$$
F_k(u_k) := f_k(p_k - g^{-1}_k (u_k)) \quad \forall k \in \{0, \dots, N-1 \},
$$
where $g^{-1}_k : [\underline{u}_k, \overline{u}_k] \mapsto [g_k^{-1}(\underline{u}_k), g_k^{-1}(\overline{u}_k)]$ models the motor's output power as a function of the battery's internal power (the inverse function exists for all $k \in \{0, \dots, N-1 \}$ \cite[Lemma 3.2.5]{east_2021}), and is defined by
$$
g_k^{-1} (u_k) : = - \frac{\beta_{1,k}}{2\beta_{2,k}} + \sqrt{-\frac{R_k u_k^2}{\beta_{2,k} V_k^2} + \frac{u_k - \beta_{0,k}}{\beta_{2,k}} + \frac{\beta_{1,k}^2}{4 \beta_{2,k}^2}}.
$$
As the control inputs are constant between update instants, the predicted internal energy of the battery, $x_k$, is given by
$$
x_{k} := x(t) - \sum_{i = 0}^{k-1} u_i \quad \forall k \in \{1, \dots, N \}.
$$
Therefore, the problem of determining the optimal battery power allocation (which implicitly determines the optimal motor and engine power)  reduces to the convex minimization
\begin{equation}\label{equation::MPC}\tag{MPC} 
\begin{aligned}
\min_{(u,x)} \ & \sum_{k \in \mathcal{P}} F_k (u_k) \\
\text{s.t.} \ & x = \mathbf{1} x(t) - \Psi u, \\
& \underline{u} \leq u \leq \overline{u},  \\
& \mathbf{1}\underline{x} \leq x \leq \mathbf{1} \overline{x},
\end{aligned}
\end{equation}
where $\underline{u}:= (\underline{u}_0, \dots, \underline{u}_{N-1})$, $\overline{u} : = (\overline{u}_0, \dots, \overline{u}_{N-1})$, $x := (x_1, \dots, x_N)$, and $\Psi$ is a lower triangular matrix of ones. Problem \eqref{equation::MPC} is convex as the objective functions $F_k$ are all convex \cite{CDC2018}, and the decision variables are only subject to affine equality and inequality constraints, which are also convex. The model predictive control input is then defined by $u^\star_0$, where $u^\star$ is the minimizing argument of \eqref{equation::MPC}. The following Lemmas are proven in \cite[\S 4.1]{east_2021}:

\begin{lemma}\label{lemma::feasible}
If \eqref{equation::MPC} is feasible, then the solution exists and is unique.
\end{lemma}

\begin{lemma}\label{lemma::feasibility}
Problem \eqref{equation::MPC} is feasible if and only if ${\mathcal{F}_k \neq \emptyset} \ \forall k \in \{1, \dots, N \}$, where $\mathcal{F}_0 : = \{x(t)\}$, and $\mathcal{F}_{k+1} : = \{ \mathcal{F}_k \oplus - [\underline{u}_k, \overline{u}_k] \} \cap [\underline{x}, \overline{x}] \ \forall k \in \{0,\dots, N-1 \} $.
\end{lemma}

\subsection{Scenario Optimization}

A limitation of the MPC energy management strategy based on the solution of \eqref{equation::MPC} is that it relies on a single, deterministic prediction of $v$ and $\theta$, which will typically be inaccurate. The ideas presented in the previous section are now extended to the stochastic case using scenario optimization.

Consider the case where $v$ and $\theta$ are drawn from sets $\mathcal{V}$ and $\Theta$ according to an associated probability. Note that the parameters $\underline{u}$ and $\overline{u}$, the functions $F_k$, and the set $\mathcal{P}$ are dependent on the realizations of $v$ and $\theta$. Therefore, the stochastic optimization problem that minimizes the expected fuel consumption across all possible realizations of $v$ and $\theta$ is
\begin{equation}\label{equation::exact_SMPC}
\begin{aligned}
\min_{(u, x)} \ & \mathbb{E}_{v \sim \mathcal{V}, \theta \sim \Theta} \left[ \sum_{\mathcal{P}(v, \theta)} F_k (u_k, v_k, \theta_k) \right] \\
\text{s.t.} \ & x = \mathbf{1} x(t) - \Psi u, \\
& \underline{u}(v, \theta) \leq u \leq \overline{u} (v, \theta) \quad \forall (v, \theta) \in \mathcal{V} \times \Theta \\
& \mathbf{1}\underline{x} \leq x \leq \mathbf{1} \overline{x}.
\end{aligned}
\end{equation}
This approach would increase the robustness of the controller as it is more likely that the resulting driver behaviour will have been explicitly considered by the MPC optimization, but still has several significant issues. Particularly, in order for the controller to minimize the expected fuel consumption when implemented in closed loop, the probability measure associated with $\mathcal{V}$ and $\Theta$ must closely match real driver behaviour, a phenomenon that is extremely challenging to model explicitly. Secondly, even in the case where an accurate representation of driver behaviour is available and the associated probability density function is `simple', problem \eqref{equation::exact_SMPC} is nonconvex and intractable in general.

Therefore, consider instead the case where $S \in \mathbb{Z}_{++}$ individual samples of $v$ and $\theta$ are obtained from the uncertainty sets, where each sampled pair $(v, \theta) \in \mathcal{V} \times \Theta$ is referred to as a \emph{scenario}. It is again assumed that each of the scenarios is sampled at $1$ Hz over the remaining duration of the journey, although as before, this approach can be extended to an arbitrary sampling interval. Additionally, it is assumed that each of the predictions are of the same length, $N$, as although the remaining duration of the journey may vary for each scenario, zeros can be appended to the shorter predictions until they equal the length of the longest prediction. Therefore, the predictions are given by $v \in \mathbb{R}^{N \times S}$ and $\theta \in \mathbb{R}^{N \times S}$, where $v_{:,s} \in \mathcal{V}$ and $\theta_{:,s} \in \Theta$ $\forall s \in \{1, \dots, S\}$. The approach detailed in Section \ref{section::convex_MPC} can then be used to obtain the optimization parameters across all of the predictions, now given by $\underline{u} \in \mathbb{R}^{N \times S}$, $\overline{u} \in \mathbb{R}^{N \times S}$, and $d \in \mathbb{R}^{N\times S}$. Additionally, the set $\mathcal{P}$ is defined for each scenario $s \in \{1, \dots, S\}$ as $\mathcal{P}_s$, and the functions $f_{k,s}(\cdot)$ and $g^{-1}_{k,s} ( \cdot )$ are defined for all elements of the prediction horizon $k \in \{0, \dots, N-1 \}$, and scenarios $s \in \{1, \dots, S \}$. The stochastic program (\ref{equation::exact_SMPC}) may then be approximated with 
\begin{equation}\tag{SMPC}\label{eqn::SMPC}
\begin{aligned}
\min_{(u,x, \tau)} \ & \frac{1}{S} \sum_{s=1}^S \bigg( \sum_{k \in \mathcal{P}_s}  f_{k,s} (p_{k,s} - g_{k,s}^{-1} (u_{k,s}) ) \bigg) , \\
\textrm{s.t.} \ 
& x_{:,s} = \mathbf{1} x(t) - \Psi u_{:,s}, \\
& \mathbf{1} \underline{x} \leq x_{:,s} \leq \mathbf{1} \overline{x}, \\
& \underline{u}_{:,s} \leq u_{:,s} \leq \overline{u}_{:,s} \\
& u_{0,s} = \tau, \\
& \forall s \in \{1, \dots, S \}.
\end{aligned}
\end{equation}
The following can be proven using minor modifications to the proofs of Lemmas \ref{lemma::feasible} and \ref{lemma::feasibility} in \cite[\S 4.1]{east_2021}:
\begin{lemma}
If \eqref{eqn::SMPC} is feasible, then the solution exists and is unique.
\end{lemma}
\begin{lemma}\label{lemma::feasibility2}
Problem \eqref{eqn::SMPC} is feasible if and only if $\mathcal{F}_{k,s} \neq \emptyset \ \forall k \in \{1, \dots, N \}$, where $\mathcal{F}_{0,s} : = \{x(t)\}$, $\mathcal{F}_{1,s} := \{ \mathcal{F}_{0,s} \oplus - \bigcap_{s=1}^S[\underline{u}_{0,s}, \overline{u}_{0,s}] \} \cap [\underline{x}, \overline{x} ] $, and $\mathcal{F}_{k+1,s} : = \{ \mathcal{F}_{k,s} \oplus - [\underline{u}_{k,s}, \overline{u}_{k,s}] \} \cap [\underline{x}, \overline{x}] \ \forall k \in \{1,\dots, {N-1} \}, \ s \in \{1, \dots, S \} $.
\end{lemma}

From a computational perspective, (\ref{eqn::SMPC}) has two significant advantages relative to (\ref{equation::exact_SMPC}). Firstly, (\ref{eqn::SMPC}) is convex and deterministic, and can therefore be readily solved using standard convex optimization algorithms. Secondly, (\ref{eqn::SMPC}) does not require the driver's behaviour to be modelled explicitly, and the scenarios can be drawn from a database obtained from previous examples of a given route (or its constituent parts, assembled from available data). An additional advantage of (\ref{eqn::SMPC}) from a control perspective is that it allows a separate control vector for each scenario (i.e. $u \in \mathbb{R}^{N \times S}$) rather than a single control vector that minimizes the objective across all predictions. This is preferable as at any given instant along the prediction horizon one scenario may predict high speed and/or acceleration whereas another may predict that the vehicle is moving slowly or is stationary (this places conflicting requirements on the predicted control input which cannot satisfy both simultaneously). The predicted control trajectories in (\ref{eqn::SMPC}) can therefore be interpreted as samples of a control probability density function, which is a more general representation than that considered in (\ref{equation::exact_SMPC}). The solution of (\ref{eqn::SMPC}) provides a meaningful MPC law since the first element of each predicted control sequence is constrained to be equal, so that the feedback law is defined by $u(t) := \tau^\star$, where $\tau^\star$ is obtained from the minimizing argument of \eqref{eqn::SMPC}. 

A limitation of (\ref{eqn::SMPC}) is that the solution only converges in the limit as $S \to \infty$ (the error term between the solution of (\ref{eqn::SMPC}) and the expected cost w.r.t.\ the true probability distribution is $\mathcal{O}(S^{-1/2})$ \cite{Shapiro1993}). This may become problematic if both a large number of scenarios and a long prediction horizon is required, as (\ref{eqn::SMPC}) has $\mathcal{O} (NS)$ decision variables. An efficient ADMM algorithm for the solution of (\ref{eqn::SMPC}) is presented in the following to address this issue. A further limitation is that the solution of (\ref{eqn::SMPC}) represents $S$ deterministic, open loop policies, so the optimal value of (\ref{eqn::SMPC}) does not necessarily converge to the resulting closed-loop cost as $S \to \infty$. Possible modifications to the proposed framework that address this optimality gap are discussed in Section \ref{subsection::results} in light of the numerical results that follow.

\subsection{Feasibility}\label{section::feasibility}

This section addresses the feasibility of (\ref{eqn::SMPC}) when implemented in closed loop. In particular, results from scenario MPC are used to obtain a bound on the number of scenarios, $S$, required to provide a given confidence that the one-step-ahead instance of (\ref{eqn::SMPC}) will be feasible with a given probability.

\begin{assumption}\label{assumption::ch6::prediction1}
The first element of the predictions of velocity and road gradient are accurate measurements of the vehicle's current velocity and the current road gradient, and are therefore equal for each scenario, i.e $v_{0,s} = v (t)$ and $\theta_{0,s} = \theta(t)$ for all $s \in \{1, \dots, S \}$.
\end{assumption}
Under Assumption \ref{assumption::ch6::prediction1}, $\underline{u}_{0,s} = \underline{u} (t)$ and $\overline{u}_{0,s} = \overline{u} (t)$ for all $s \in \{1, \dots, S \}$ (where $\underline{u}(t) \in \mathbb{R}$ and $\overline{u}(t) \in \mathbb{R}$ are the lower and upper limits on battery power at time $t$).
\begin{assumption}\label{assumption::power1}
$\underline{u}(t) \leq \overline{u}(t)$ for all $t \in \{0, \dots, T \}$.
\end{assumption}
Assumption \ref{assumption::power1} can reasonably be expected to hold in general using the method outlined in \cite{CDC2018}, and must hold in practice (the battery must do \emph{something}). Consequently, Lemma \ref{lemma::feasibility2} implies that a given instance of \eqref{eqn::SMPC} is feasible if and only if
\begin{gather}
x(t) - \overline{u}(t) \leq \overline{x}, \label{gather::1} \\
x(t) - \underline{u}(t) \geq \underline{x}, \label{gather::2}  \\
x(t) - \overline{u}(t) - \sum_{i=1}^{k} \overline{u}_{i,s} \leq \overline{x},  \label{gather::3}  \\
x(t) - \underline{u}(t) - \sum_{i=1}^{k} \underline{u}_{i,s} \geq \underline{x}, \label{gather::4} 
\end{gather}
for all $k \in \{1, \dots, N - 1 \}$ and $s \in \{1, \dots, S \}$. The physical interpretation of the conditions \eqref{gather::1} and \eqref{gather::3} being violated is that so much regenerative braking is used/expected that the powertrain will be forced to overcharge the battery. It can reasonably be assumed that any given vehicle will be designed so that all of the required braking power can be delivered by the mechanical brakes, so it therefore follows that the power allocated to the electric system can always be modified to ensure that \eqref{gather::1} and \eqref{gather::3} hold. However, driveability issues may be introduced if the relative fraction of braking power required from the motor and mechanical brakes changes significantly from one MPC update to the next, in order to ensure that \eqref{gather::1} holds. The physical interpretation of the conditions \eqref{gather::2} and \eqref{gather::4} being violated is that so much drive power is demanded/expected from the powertrain that the battery will be forced to over-discharge. In this case, either the engine needs to be switched on to help meet the driver's power demand (which can also introduce driveability issues if this action is taken to ensure \eqref{gather::2}), or, if the engine is already on, the power delivered by the powertrain will need to be reduced below that demanded by the driver (which could be dangerous if it occurs, for example, during an overtaking maneuver).
\begin{assumption}\label{assumption::scenario1}
The scenarios are generated so that \eqref{gather::3} and \eqref{gather::4} always hold.
\end{assumption}
Assumption \ref{assumption::scenario1} reflects the fact that the above considerations can be taken into account when determining the engine switching and braking control used for the predictions. Consequently, problem (\ref{eqn::SMPC}) is feasible at time $t$ if and only if \eqref{gather::1} and \eqref{gather::2} hold, and feasible at time $t+1$ if and only if   
\begin{equation}\label{eqn::ch6::feasability_ahead}
x(t) - \tau^\star(t) - \overline{u}(t+1) \leq \overline{x} \quad \textrm{and} \quad x(t) - \tau^\star(t) - \underline{u}(t+1) \geq \underline{x},
\end{equation}
where $u^\star(t)$ and $\tau^\star (t)$ are obtained from the solution of (\ref{eqn::SMPC})  at time $t$. Given the above physical interpretations of these constraints, the concept of feasibility here is synonymous with \emph{the powertrain operating as intended at the current time}. 

Now, note that problem (\ref{eqn::SMPC})  is equivalent to
\begin{subequations}\label{eqn::ch6::SMPC}
\begin{align}
\min_{(u,x, \tau) \in \mathcal{C}} \ & \frac{1}{S} \sum_{s=1}^S \bigg( \sum_{k \in \mathcal{P}_s}  f_{k,s} (p_{k,s} - g_{k,s}^{-1} (u_{k,s}) ) \bigg) , \\
\textrm{s.t.} \ 
& x(t) - \tau \in [\underline{x} + \underline{u}_{1,s}, \overline{x} + \overline{u}_{1,s}] \label{subeqn::ch6::feasiility} \\
& \forall s \in \{1, \dots, S \}, \notag
\end{align}
\end{subequations}
(where $\mathcal{C}$ is a convex constraint set), which is the form of problem considered in \cite{Schildbach2013}. To apply the results of~\cite{Schildbach2013}, assume that (\ref{eqn::SMPC}) is feasible at time $t$, and define the violation probability for the one-step-ahead feasibility condition (\ref{eqn::ch6::feasability_ahead}) as
$$
V(\tau^\star(t)) := \mathbb{P} \{ x(t) - \tau^\star(t) \notin [\underline{x} + \underline{u}(t+1), \overline{x} +  \overline{u}(t+1)] \},
$$
where $\underline{u}(t+1)$ and $\overline{u}(t+1)$ can be interpreted as an unseen sample of the random variables $\underline{u}_{1,s}$ and $\overline{u}_{1,s}$ in (\ref{eqn::SMPC}). The support rank of constraint (\ref{subeqn::ch6::feasiility}) is one\footnote{ Note that, whist the set $\mathcal{C}$ depends on the $S$ random samples of $\mathcal{V}$ and $\Theta$ and is therefore `stochastic', the problem of interest here is whether the solution of (\ref{eqn::SMPC}) remains feasible with respect to an additional randomly sampled constraint of the form of (\ref{subeqn::ch6::feasiility}) \emph{only}, and not every constraint from (\ref{eqn::SMPC}) (i.e. $\mathcal{C}$ does not change when the new scenario is sampled). Therefore, $\mathcal{C}$ can be considered deterministic, and the relevant support rank is that of constraint (\ref{subeqn::ch6::feasiility}).}, since the constraint fixes a single degree of freedom of the decision variable $(u, x, \tau)$ (see Definition 3.6 in \cite{Schildbach2013}). Therefore, the confidence bound
\begin{equation}\label{eqn::ch6::confidence_bound}
\mathbb{P} \{V(\tau^\star (t)) > \hat{\epsilon} \} \leq B(\epsilon; S, 0) =  (1 - \hat{\epsilon} )^S
\end{equation}
 on the violation (i.e. one-step-ahead infeasibility) probability of \eqref{eqn::SMPC} immediately follows from Theorem 4.1 in \cite{Schildbach2013}. More generally, if $R$ samples are discarded (using any criteria, e.g. a greedy strategy), then the right-hand-side of (\ref{eqn::ch6::confidence_bound}) can be replaced with $B(\hat{\epsilon};S,R)$, where $B$ is the beta distribution.

Let the confidence that the violation probability is no greater than $\hat{\epsilon}$ be $\hat{\beta}$, so that
\begin{equation}\label{eqn::ch6::confidence_bound2}
\mathbb{P} \{ V ( \tau^\star(t)) \leq \hat{\epsilon} \} = \hat{\beta} \geq 1 - (1 - \hat{\epsilon} )^S,
\end{equation}
then the lower bound on the number of scenarios required to meet a given confidence of violation probability is 
\begin{equation}\label{eqn::ch6::number_of_scenarios}
S \geq \frac{\log{(1 - \hat{\beta})}}{\log{(1 - \hat{\epsilon } ) }}.
\end{equation}
This bound is illustrated in Figure \ref{figure::violation_probability}. One-step-ahead feasibility (and by induction, recursive feasibility) is only guaranteed with complete certainty in the limit as $\mathcal{S} \to \infty$. However, note that (\ref{eqn::ch6::confidence_bound2}) is a \emph{lower bound} on the confidence of constraint violation probability, which will only be tight when the state-of-charge of the battery is `close' to its limits (i.e. one-step-ahead feasibility is guaranteed when the battery's state-of-charge is `far' from its limits, as it will not be possible to exceed the limits within two sampling intervals). It is not common for the state-of-charge of the battery to approach its upper or lower bounds more than once or twice during a given journey (see for example \cite[\S 5]{Sun}, \cite[\S 5]{Stockar}, \cite[\S 5]{Hu2013}, and the results that follow), so high confidence of feasibility throughout a given journey may be obtained with a modest number of scenarios. For example, a bound on constraint violation of $\hat{\epsilon} =  0.1$ can be achieved with confidence of $\hat{\beta} \geq 0.9$ using just $S = 22$ scenarios. This can be improved to $V ( \tau^\star(t)) \leq 0.01$ with a confidence of $\hat{\beta} \geq 0.99$ by using $S = 459$ scenarios. 

\usepgfplotslibrary{fillbetween}
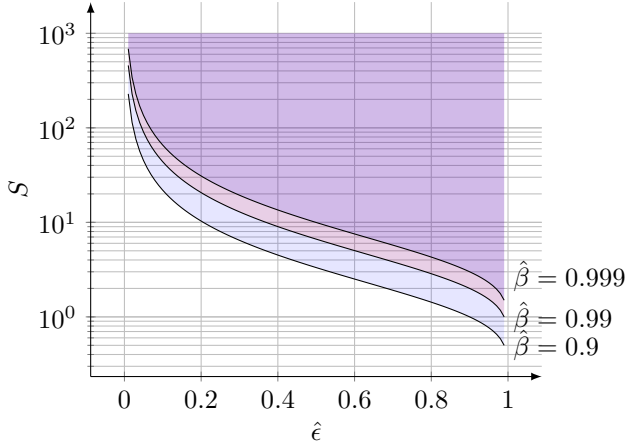
\begin{figure}
\begin{center}
\begin{tikzpicture}[>=latex, scale=0.875]
\begin{axis}[xlabel=$\hat{\epsilon}$,ylabel=$S$,axis x line=bottom, axis y line = left, axis line style={->}, enlargelimits, ymode=log, clip=false, grid=both]
\addplot[domain=0.01:0.99, samples=100, name path=A, opacity = 0] {1000};
\addplot[domain=0.01:0.99, samples=100, name path=B] {ln(1 - 0.9) / ln(1 - x)} node[right, pos=1] {$\hat{\beta} = 0.9$};
\addplot[domain=0.01:0.99, samples=100, name path=C] {ln(1 - 0.99) / ln(1 - x)} node[right, pos=1] {$\hat{\beta} = 0.99$};
\addplot[domain=0.01:0.99, samples=100, name path=D] {ln(1 - 0.999) / ln(1 - x)} node[above right, pos=1] {$\hat{\beta} = 0.999$};
\addplot [blue, opacity=0.1] fill between [of=A and B];
\addplot [red, opacity=0.1] fill between [of=A and C];
\addplot [blue, opacity=0.1] fill between [of=A and D];
\end{axis}
\end{tikzpicture}
\caption{Illustration of the values of $S$ that satisfy (\ref{eqn::ch6::number_of_scenarios}) for $\hat{\beta} \in \{0.9, 0.99, 0.999 \}$ and $\hat{\epsilon} \in [0.01, 0.99]$. \label{figure::violation_probability}}
\end{center}
\end{figure}

\section{Optimization Algorithm}
\label{section::optimization}

 It was shown in \cite{EastTCST} that ADMM can be used to solve the deterministic case \eqref{equation::MPC} in a fraction of a second for prediction horizons exceeding 1000 time steps, and that it is generally orders of magnitude faster than the general purpose convex optimization software CVX \cite{CVX} with SDPT3 \cite{SDPT3}. Therefore, in this section we extend the ADMM algorithm in \cite{EastTCST} to the scenario optimization \eqref{eqn::SMPC}.

\subsubsection{Algorithm}

Problem (\ref{eqn::SMPC}) is equivalent to the equality constrained problem
\begin{equation}\label{eqn::ch6::opt_equality_constrained}
\begin{aligned}
\min_{(u, x, \tau)} \ & \frac{1}{S} \sum_{s=1}^S \Bigg( \sum_{k \in \mathcal{P}_s} \bigg( f_{k,s} (p_{k,s} - g_{k,s}^{-1} (u_{k,s}) ) \bigg) \\
& \hspace{1cm} + \sum_{k=0}^{N-1} \left( \mathcal{I}_{\U_{k,s}}(u_{k,s}) + \mathcal{I}_{\X} (x_{k+1,s}) \right) \Bigg), \\
\textrm{s.t.} \ 
& u = \zeta, \quad x = \Phi x(t) - \Psi \zeta, \quad \mathbf{1} \tau = u_{0,:}^\top, \\
\end{aligned}
\end{equation}
which is in turn equivalent to
\begin{equation}\label{eqn::ch6::opt_canonical}
\min_{\tilde{u}, \tilde{x}} \ \tilde{f} ( \tilde{u} )  \quad \textrm{s.t } A \tilde{u} + B \tilde{x} = c,
\end{equation}
where
\begin{gather*}
\tilde{u} := (\vec (u), \vec (x)), \quad \tilde{x}: = (\tau, \vec (\zeta) ),  \\
\tilde{f}(\tilde{u}) : =  \frac{1}{S} \sum_{s=1}^S \Bigg( \sum_{k \in \mathcal{P}_s} \bigg( f_{k,s} (p_{k,s} - g_{k,s}^{-1} (u_{k,s}) ) \bigg) \\ 
\hspace{3cm} + \sum_{k=0}^{N-1} \left( \mathcal{I}_{\U_{k,s}}(u_{k,s}) + \mathcal{I}_{\X} (x_{k+1,s}) \right) \Bigg), \\
A:= \left[ \begin{smallmatrix}
I & O \\
O & I \\
-\left( \mathbf{1} \otimes [1,0,\dots,0] \right) & O
\end{smallmatrix} \right] , \quad B:= \biggl[ \begin{smallmatrix}
O & -I \\
O  & (I \otimes \Psi ) \\
 \mathbf{1} & O
\end{smallmatrix} \biggr], \\
c:= (\mathbf{0}, \mathbf{1}x(t), \mathbf{0} ).
\end{gather*}
Now define the augmented Lagrangian function $\mathcal{L}_\rho :  \mathbb{R}^{2NS} \times \mathbb{R}^{1+NS} \times \mathbb{R}^{(2N + 1 )S} \mapsto \mathbb{R}$ as
\begin{multline*}
\mathcal{L}_\rho (\tilde{u}, \tilde{x}, \lambda ) := \frac{1}{S} \sum_{s=1}^S \Bigg( \sum_{k \in \mathcal{P}_s} \bigg( f_{k,s} (p_{k,s} - g_{k,s}^{-1} (u_{k,s}) ) \bigg) + \\
\sum_{k=0}^{N-1} \left( \mathcal{I}_{\U_{k,s}}(u_{k,s}) + \mathcal{I}_{\X} (x_{k+1,s}) \right) + \frac{\rho_1}{2} \| u_{:,s} - \zeta_{:,s} + \nu_{:,s} \|^2  \\ 
+ \frac{\rho_2}{2} \| x_{:,s} + \Psi \zeta_{:,s} - \mathbf{1} x(t) + \psi_{:,s} \|^2 + \frac{\rho_3}{2} (\tau - u_{0,s} + \phi_{s} )^2 \Bigg),
\end{multline*}
where $\lambda:= (\vec (\nu), \vec (\psi), \phi)$ is a vector of Lagrange multipliers and $\nu \in \mathbb{R}^{N\times S}$, $\psi \in \mathbb{R}^{N \times S}$, $\phi \in \mathbb{R}^S$. Problem (\ref{eqn::ch6::opt_canonical}) is the canonical form of the optimization problem considered in \cite{ADMM}, where it is shown that the ADMM iteration 
\begin{subequations}\label{eqn::ADMM}
\begin{align}
\tilde{u}^{(i+1)} &: = \argmin_{\tilde{u}} \mathcal{L}_\rho (\tilde{u}, \tilde{x}^{(i)}, \lambda^{(i)} ) \label{subeq::1} \\
\tilde{x}^{(i+1)} &: = \argmin_{\tilde{x}} \mathcal{L}_\rho (\tilde{u}^{(i+1)}, \tilde{x}, \lambda^{(i)} ) \label{subeq::2} \\
\lambda^{(i+1)} & := \lambda^{(i)} + A \tilde{u}^{(i+1)} + B \tilde{x}^{(i+1)} - c \label{subeq::3}
\end{align}
\end{subequations}
converges\footnote{The proof is provided for a single, scalar value of $\rho$, but trivially extends to the case here, where $\rho$ can be represented by a diagonal matrix with positive elements. Also, note that \eqref{eqn::SMPC} has only affine inequality constraints, so feasibility is a sufficient condition for strong duality \cite[Proposition 5.3.1]{bertsekas_convex_opt_theory}.} to the (unique, in this case) solution. The optimality and convergence of a given iterate $(\tilde{u}^{(i+1)}, \tilde{x}^{(i+1)}, \lambda^{(i + 1)})$ is indicated by the sizes of the residual variables defined by
\begin{align}\label{eqn::residuals}
r\powip & := \begin{bmatrix} \vec \left( u\powip \right) - \vec \left( \zeta\powip \right) \\ \vec x - \mathbf{1} x(t) + (\Psi \zeta_{:,1}\powip , \dots, \Psi \zeta_{:,S}\powip ) \\
 \mathbf{1} \tau\powip - (u_{0,:} \powip)^\top \end{bmatrix} \\
s\powip & := \begin{bmatrix}
\rho_1  \vec (\zeta \powi - \zeta \powip ) - M (\tau\powi - \tau \powip ) \\
\rho_2 \Big( \Psi ( \zeta_{:,1}\powi - \zeta_{:,1}\powip), \dots, \Psi ( \zeta_{:,S}\powi - \zeta_{:,S}\powip ) \Big)
\end{bmatrix},
\end{align}
where $M:=\rho_3 ( \mathbf{1}_S \otimes [1, 0, \dots, 0]^\top )$. The algorithm is terminated at the first iteration in which the criterion 
$$
\max \{ \| r^{(i+1)} \|, \| s^{(i+1)} \| \} \leq \epsilon
$$
is met, where $\epsilon \in \mathbb{R}_{++}$ is a pre-determined convergence threshold. The iteration \eqref{eqn::ADMM} can be made to terminate at a point that can be made arbitrarily close to the unique minimizing argument of (\ref{eqn::SMPC}) by setting the termination threshold, $\epsilon$, arbitrarily close to zero.

\subsection{Variable Updates \& Complexity}\label{section::ch6::complexity}

The $\tilde{u}$ ADMM update \eqref{subeq::1} is equivalent to 
\begin{align*}
u_{0,s}\powip & = \argmin_{u_{0,s}} \Big(  f_{0,s} (p_{0,s} - g_{0,s}^{-1} (u_{0,s}) ) + \mathcal{I}_{\U_{0,s}}(u_{0,s}) \\
& \hspace{-0.15cm} + \frac{\rho_1}{2} (u_{0,s} - \zeta_{0,s}\powi + \nu_{0,s}\powi )^2 + \frac{\rho_3}{2} (\tau\powi - u_{0,s} + \phi_s\powi)^2 \Big) \\
& \forall s \in \{1, \dots, S \} \\
u_{k,s}\powip & = \argmin_{u_{k,s}} \Big(  f_{k,s} (p_{k,s} - g_{k,s}^{-1} (u_{k,s}) ) + \mathcal{I}_{\U_{k,s}}(u_{k,s}) \\
& \hspace{-0.15cm}  + \frac{\rho_1}{2} (u_{k,s} - \zeta_{k,s}\powi + \nu_{k,s}\powi )^2 \Big) \quad  \\
& \forall s \in \{0, \dots, S \}, \ k \in \mathcal{P}_s \setminus \{0\} \\
u_{k,s}\powip & = \underline{u}_{k,s} \quad  \forall s \in \{0, \dots, S \}, k \notin \mathcal{P}_s  \\
x_{:,s}\powip & = \Pi_{\X^N} \left[ \mathbf{1} x(t) -  \Psi \zeta_{:,s}\powi - \psi_{:,s}\powi \right],
\end{align*}
where $\mathcal{X} := \{x \in \mathbb{R} : \underline{x} \leq x \leq \overline{x} \} $, and $\Pi_{\mathcal{X}^N} (x):= \argmin_{\hat{x} \in \mathcal{X}^N} \| x - \hat{x} \|$. The update $u_{k,s}\powip$ is a 1-dimensional convex inequality constrained optimization problem that can be solved using Newton's method for all $s \in \{1, \dots, S\}$ and $k \in \mathcal{P}_s$. This update is separable w.r.t.\ both the scenarios and the prediction horizon, and so the computational requirement is constant as both are increased if the updates are performed in parallel (so long as $\sum_{s=1}^S |\mathcal{P}_s|$ threads are available). Otherwise, the computational complexity is $\mathcal{O}(NS)$ if the updates are performed sequentially.
The computation of $x\powip$ is $\mathcal{O}(N)$ as multiplication by $\Psi $ is a cumulative sum and  projection $\Pi_{\mathcal{X}^N}[\cdot]$ is performed element-wise. No additional memory is required for either the $u\powip$ or the $x \powip$ update.

The $\tilde{x}$ ADMM update (\ref{subeq::2}) is equivalent to
\begin{align*}
\tau\powip & = \frac{1}{S} \sum_{s=1}^S (u_{0,s}\powip + \phi_s\powi ) \\
\zeta_s\powip & = (\rho_1 I + \rho_2 \Psi^\top \Psi )^{-1}  \\
& \hspace{0.9cm} \Big( \rho_1 (u_{:,s} + \nu_{:,s}) - \rho_2 \Psi^\top ( x_{:,s} - \mathbf{1} x(t) + \psi_{:,s} ) \Big).
\end{align*}
The computation of $\tau\powip$ scales linearly with $S$ regardless of whether parallel processing is used, but since this is a scalar sum it is typically computationally inexpensive. The computation required to solve the linear system of equations $(\rho_1 I + \rho_2 \Psi^\top \Psi ) \mathbf{x} = \mathbf{b}$ (with $\mathbf{x} \in \mathbb{R}^N$, $\mathbf{b} \in \mathbb{R}^N$) and evaluate the matrix-vector product $\Psi^\top \mathbf{x}$ scales linearly with $N$ (the proof of the former is available in \cite[Appendix D]{EastTVT}, while multiplication by $\Psi^\top$ is a cumulative sum). Therefore the computation of $\zeta^{(i+1)}$ scales as $\mathcal{O}(N)$ if the $S$ updates are performed in parallel, or $\mathcal{O}(NS)$ if they are performed sequentially. The Cholesky factors required for the solution of the linear system of equations (see  \cite[Appendix D]{EastTVT} for details) are the same for each scenario so only need storing once. Therefore, the memory cost of the update is $\mathcal{O}(N)$. 

The computation of residual updates \eqref{eqn::residuals} scales linearly with $N$ and $S$ (the only non-trivial computations are multiplication by $\Psi$, which is equivalent to a cumulative sum, and multiplication by $(\mathbf{1}_S \otimes [1,0, \dots, 0]^\top)$, which is a mathematical representation of a simple operation that copies the value of $(\tau^{(i)} - \tau^{(i+1)})$ into entries of a sparse column vector). Additionally, $\mathcal{O}(NS)$ memory storage is required to store the variables $\tau \powi$ and $\zeta\powi$ for the computation of $s\powip$ at the following iteration. The $\lambda$ update is equivalent to $\lambda\powip = \lambda\powi + r\powip$.

\section{Numerical Experiments}
\label{section::numerical_experiments}

The performance of the proposed SMPC was investigated using numerical simulations of a PHEV completing 49 instances of a 13 km route, as illustrated in Figure \ref{fig:driver_behaviour}. The key property of this data is that is represents multiple trips of the same route, and can therefore be readily implemented in the scenario MPC framework. Additionally, it was obtained using four different drivers, in weather conditions including `dry', `wet', and `foggy', which therefore implies that an energy management strategy that performs well across all journeys is robust to variations in driving style and weather. Detailed descriptions of the traffic conditions are unavailable, but start-stop artefacts are absent from the data, suggesting that the traffic was generally light. The route is a 6.5 km return trip, and includes a junction (visible from the speed drop at $\sim$5 km and returning at $\sim$9 km) and a roundabout at the midway point (visible from the speed drop at 6.5 km). The first half of the route is predominantly uphill, and therefore the second half is predominantly downhill. The PHEV was modelled for simulation as a 1800 kg passenger vehicle with a 100 kW engine, 50 kW motor, and a 21.5 Ah battery, which is representative of the vehicle used to obtain the data. For each journey, two controllers were used: 

\begin{figure}
    \begin{center}
    \input{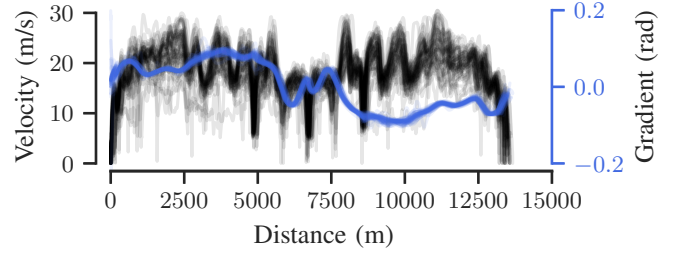}
    \end{center}
    \caption{Velocity and road gradient against longitudinal distance for each of the 49 journeys.}
    \label{fig:driver_behaviour}
\end{figure}

\begin{enumerate}
\item A \textbf{Nominal MPC} law defined by $u(t) := u^\star_0$, where $u^\star_0$ is the first element of the minimizing argument of \eqref{equation::MPC}. At each control variable update instance, the MPC used a single scenario that was a completely accurate representation of driver behaviour during the entire remaining portion of the journey, and was solved using the ADMM algorithm from \cite{EastTCST}. 
\item A \textbf{Scenario MPC}
law defined by $u(t) := \tau^\star$ where $\tau^\star$ was obtained from the solution of (\ref{eqn::SMPC}) using the ADMM algorithm described in Section \ref{section::optimization}. At each control variable update instance, the MPC used scenarios generated from the 48 examples of driver behaviour \emph{other} than the one being used for simulation; the full method is described in the following section.
\end{enumerate}
For both controllers, the gear ratio was selected using a simple vehicle speed-based heuristic, the engine was always switched on if the rotational speed of the engine was greater than its minimum speed (and otherwise switched off), and the mechanical braking power was set as a fixed fraction of braking power. The choice of these heuristics is not important for the purposes of this paper; only that the same heuristics are used for both MPC approaches. For each journey, the state of charge of the battery was initialized at 60\% and constrained between 40\% and 100\%. These simulation parameters were chosen to ensure that the lower state-of-charge constraint was active during the simulations, as the solution of (\ref{eqn::SMPC}) becomes significantly harder under this condition. The parameters $\rho_1$ and $\rho_2$ in the ADMM algorithm were set at $\rho_1 = 2.34 \times 10^{-4}$ and $\rho_2 = 8.86 \times 10^{-9}$ (as in \cite{CDC2018}), and $\rho_3$ was set as equal to $\rho_1$ (as they both correspond to constraints on power). The termination criterion was set at $\epsilon = 0.1$ (this value was found to be sufficiently low to avoid the constraint satisfaction issue discussed in \cite[\S 5-C]{EastTCST}). The simulations were programmed in Matlab and
run on a 2.60GHz Intel i7-9750H CPU.

\subsection{Scenario Generation}

A database was available to the controller where arrays of velocity, road gradient, and time data, indexed by the longitudinal position within the route, $l$, and sampled at a resolution of 1m, was available for all of the 49 journeys used for simulation. Before the start of each simulated journey, the position indexed data for $S = 48$ journeys (i.e. every journey \emph{except} the one representing the actual journey to be simulated) was retrieved. The data used for simulation was indexed by time, and at each sample time $t \in \{0 , \dots , T - 1 \}$ during a simulated journey, the scenarios were generated using the following steps:
\begin{itemize}
\item The longitudinal position of the vehicle was obtained from $l(t) = \sum_{i=0}^{t-1} v(i)$. 
\item For each of the $S = 48$ journeys that had been retrieved from the database, a scenario was generated by:
\begin{itemize}
\item Re-sampling the position-indexed velocity ($\hat{v}$) and road gradient ($\hat{\theta}$) data for the remaining simulation distance in terms of time at a rate of 1 Hz.
\item Calculating the scenario road gradient and velocity values using the `blending' coefficient $e^{-0.25k}$ and\footnote{The acceleration values used in the longitudinal power model (\ref{equation::longitudinal_model}) were obtained from $\dot{v}_{k,s} := \dot{v}(t) e^{-0.25k} + \dot{\hat{v}}_k ( 1 -  e^{-0.25k} )$ for all $k \in \{0, \dots, N-1\}$, where $\dot{v}(t)$ is the vehicle's simulated acceleration at time $t$ and $\dot{\hat{v}}_k$ is obtained by numerically differentiating $\hat{v}$. The values of $v_{:,s}$ were not numerically differentiated directly as additional acceleration may be introduced by (\ref{eqn::ch6::blending}) which would distort the predicted power demand.} 
\begin{equation}\label{eqn::ch6::blending}
\begin{bmatrix}
\theta_{k,s} \\
v_{k,s} \\
\end{bmatrix} : = \begin{bmatrix}
\theta (t) \\
v (t) \\
\end{bmatrix} e^{- 0.25 k} + \hspace{-2pt} \begin{bmatrix}
\hat{\theta}_{k,s} \\
\hat v_{k,s} \\
\end{bmatrix} (1 - e^{-0.25 k} ),
\end{equation}
$\forall k \in \{ 0, \dots, N-1 \}$.
\end{itemize}
\end{itemize}

\begin{figure}
    \begin{center}
    \input{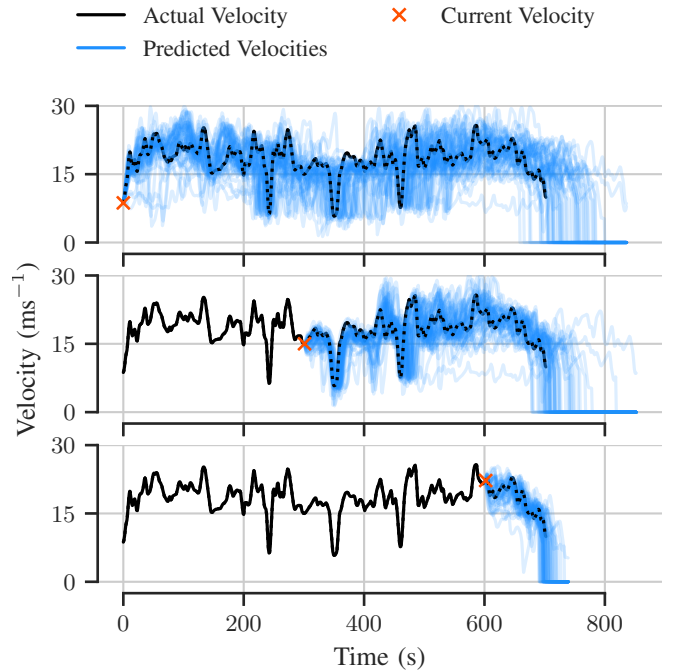}
    \end{center}
    \caption{An example of the velocity predictions (blue lines) generated at three instants during an individual journey (shown in black).}
    \label{fig:predictions}
\end{figure}

The principle of this approach is that the predictions are made \emph{entirely} using other data `previously' recorded from the route, and that the controller has no additional information about the current journey. The `blending' process (\ref{eqn::ch6::blending}) was included as it was found that the near horizon predictions (i.e. $k \in \{0, \dots, 5 \}$) were otherwise highly inaccurate due to the uncertainty in the velocity at any given point in the journey, as shown in Figure  \ref{fig:driver_behaviour}. An example of the predictions generated using the above process is illustrated in Figure \ref{fig:predictions}. It can be observed in Figure \ref{fig:driver_behaviour} that whilst there is uncertainty in the velocity at any given point on the route, there is also some discernible structure: for example, there is a roundabout at the midpoint of the journey, and the velocity in all examples drops at this point. It is interesting to note that the long-horizon predictions presented in the top plot of Figure \ref{fig:predictions} do not retain this structure, and the predictions at the later time steps do not correlate with the simulated velocities at all. It can, however, also be seen that as the vehicle approaches certain points in the journey (corresponding to locations at which the data has reduced variability), the mid-horizon predictions begin to accurately approximate the simulated velocity. For example, it can be seen in the middle plot that the first $\sim$100s of the predictions accurately match the simulated profile.

\subsection{Results}\label{subsection::results}

\begin{figure*}
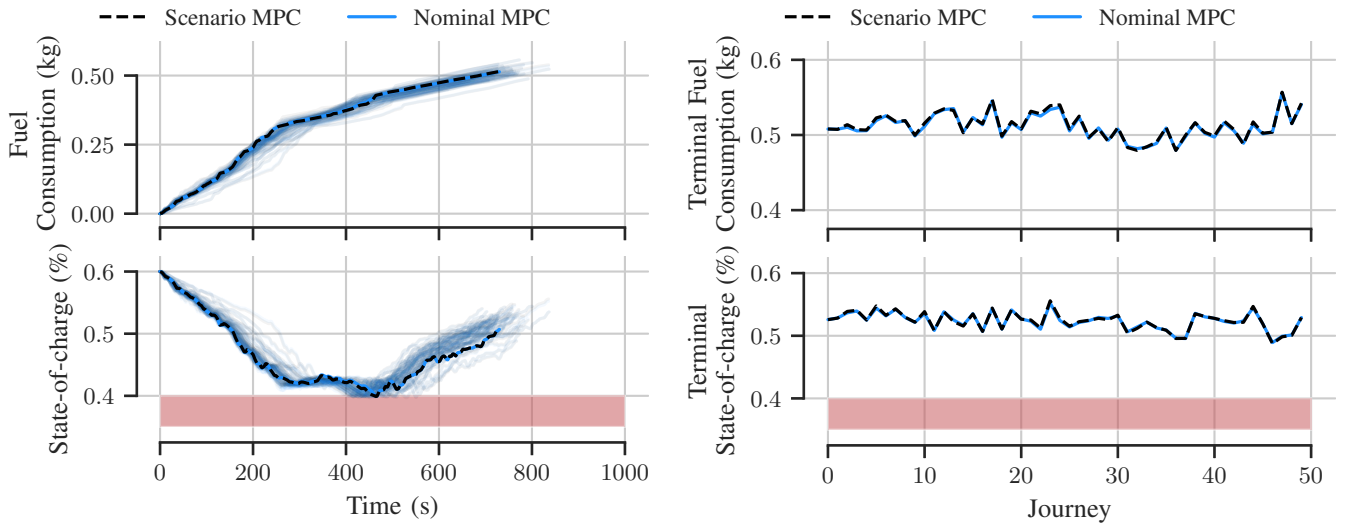

    \centering
    \begin{subfigure}[t]{0.45\textwidth}
	    \input{results_time.pgf}
	    \caption{Fuel consumption and state-of-charge trajectories with respect to time. \label{fig:results_time_1}} 
	\end{subfigure}
	\hspace{20pt}
	\begin{subfigure}[t]{0.45\textwidth}
    \centering
	    \input{results_time2.pgf}
    \caption{Fuel consumption and state-of-charge at the end of each journey. \label{fig:results_time_2}}
	\end{subfigure}
	\caption{Closed-loop fuel consumption and state-of-charge results using nominal and scenario MPC for all simulated journeys. \label{fig:results_time}}
\end{figure*}

Figure \ref{fig:results_time} shows the closed-loop simulation results using both nominal MPC controller and scenario MPC. A useful artefact of the test data is that the vehicle is travelling downhill (and therefore regenerating) for the second half of the journey, so each journey terminates with the same state-of-charge when using both controllers, and the comparative optimality can be determined by comparing fuel consumption directly. The results are remarkable: the scenario MPC state-of-charge trajectories and fuel consumption closely match those obtained using nominal MPC, and the total fuel consumption is equal (to within meaningful precision) across all 49 journeys. The suboptimality introduced to the nominal MPC controller by the model approximations in the engine, motor, and battery are also negligible (see the results obtained using dynamic programming for the same dataset in \cite[\S 4.3.2]{east_2021}), so the scenario MPC also obtains an extremely close approximation of the \emph{globally} optimal solution. This is made more remarkable by the fact that the driver behaviour is collected from 4 different drivers, which implies that modelling an individual driver's behaviour may not be necessary to achieve optimal performance.

\begin{figure*}
    \centering
    \begin{subfigure}[c]{0.45\textwidth}
    	\centering
	    \input{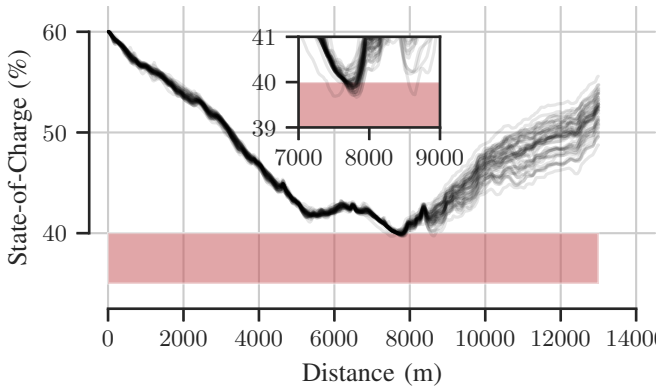}
    	\caption{Results from Figure \ref{fig:results_time_1} (using $S=48$ scenarios).}
    	\label{fig:results_distance_1}
	\end{subfigure}
	\hspace{20pt}
	\begin{subfigure}[c]{0.45\textwidth}
    	\centering
	    \input{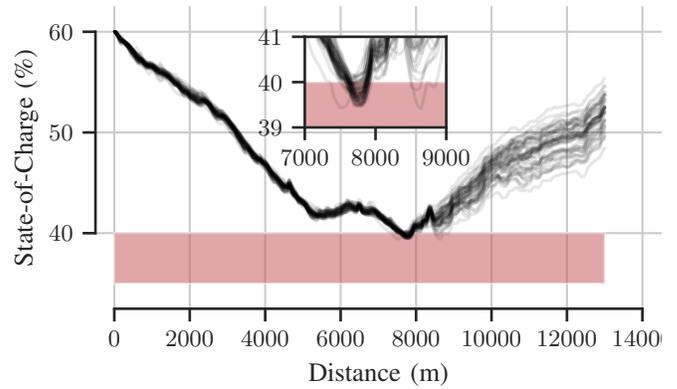}
    	\caption{Simulations from Figure \ref{fig:results_time_1} repeated using $S=5$ scenarios.}
    	\label{fig:results_distance_2}
	\end{subfigure}
	\caption{Closed-loop state-of-charge trajectories obtained using scenario MPC re-sampled w.r.t. completed distance.}
\end{figure*}

An explanation for this performance is that, although the optimal state-of-charge trajectories vary with respect to time, they match extremely closely with respect to distance. This is illustrated in Figure \ref{fig:results_distance_1}, where the closed-loop state-of-charge trajectories from Figure \ref{fig:results_time_1} are re-sampled with respect to the completed distance of the journey, and can be seen to match closely during the first half of the journey when the car is driving uphill and drive power is required from the engine. These results suggest that energy management solutions that consider driver uncertainty may be improved by considering \emph{completed distance} as the independent variable in the optimization, instead of the \emph{time}-based formulation that is overwhelmingly common in the literature. For the approach presented here, it may be possible to reformulate \eqref{eqn::SMPC} in this way with a single control trajectory across the prediction horizon, which would therefore ensure that the open loop predicted cost would converge to the closed-loop as $S \to \infty$, and reduce the computational cost of the solution. However, the trajectories in Figure \ref{fig:results_distance_1} diverge in the second half of the journey. This is because the vehicle is descending a hill and is either braking or coasting, so the energy that can be regenerated depends on losses such as aerodynamic drag. It is possible that a route that is not as evenly partitioned into an ascending phase and a descending phase may cause the scenario controller to behave more poorly. A systematic investigation into these possibilities and issues is an interesting direction for future work.

The results in Figure \ref{fig:results_distance_1}  show that the lower state-of-charge constraint is violated by up to $\sim 0.2 \%$ (if the constraints were violated then the simulation was allowed to continue, and the behaviour of the driver was not modified as discussed in Section \ref{section::feasibility}).  The bound (\ref{eqn::ch6::confidence_bound2}) provides a confidence of $\hat{\beta} \geq 0.92$ for a violation probability of $\hat{\epsilon} = 0.05$ when using $S = 48$ scenarios, so a portion of this constraint violation can be attributed to an insufficient number of predictions. This is supported by the results in Figure \ref{fig:results_distance_2}, which show the results from Figure \ref{fig:results_distance_1} repeated with only $S=5$ scenarios, and for which the constraint violation increases up to a maximum of $\sim 0.5 \%$ (more scenarios could not be tested due to a lack of data). It is worth noting, however, that the bound (\ref{eqn::ch6::confidence_bound2}) is dependent on the scenarios being sampled from the uncertainty set $\mathcal{V} \times \Theta$, and the one-step-ahead predictions obtained using the blending process (\ref{eqn::ch6::blending}) are not empirical. Therefore, it may also be possible to reduce the rate and severity of constraint violations by sampling the one-step-ahead predictions alone from a separate distribution, possibly conditioned on the current road gradient, and vehicle velocity and acceleration.  Nonetheless, the observed constraint violation is likely to be lower than the uncertainty in the battery state estimate, so may not become an issue in practice. 

\begin{figure}
    \centering
	    \input{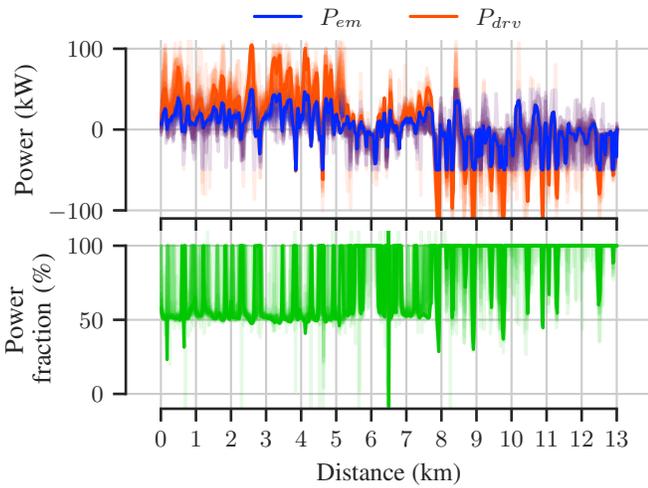}
    \caption{Control decision variables for each of the journeys illustrated in Figure \ref{fig:results_time} using scenario MPC, with the same particular journey highlighted. The top axes show the power demanded by the driver and the power delivered by the motor, and the bottom axes show the power delivered by the motor as a fraction of the driver's demanded power.}
    \label{fig:results_distance_control}
\end{figure}

Figure \ref{fig:results_distance_control} illustrates the power-split control actions for each the journeys when using scenario MPC, and includes the power demanded by the driver, $P_{drv}$, the power delivered by the motor, $P_{em}$, and the fraction of power delivered by the motor relative to the total demand power, all sampled w.r.t. distance (these results are not illustrated w.r.t. time as they become indecipherable). Again, these results show a strong correlation between the control actions determined for each journey at any given position along the route. An interesting additional artefact is that during the first half of the journey (where the vehicle is ascending the hill), the power fraction delivered by the motor generally switches between 100\% and $\sim$50\%. The portions where the motor is delivering 100\% of the powertrain's output generally correspond to sections of the journey where the power demand drops to a very low level or braking occurs, i.e. the regions where the motor is \textit{required} to deliver/accept all of the driver's demand power. These results appear to support the rationale of equivalence factor-based approaches to power-split optimization \cite{Sciarretta}, although further analysis is required to determine if this behaviour is observed more generally.

\begin{figure}
    \centering
    \input{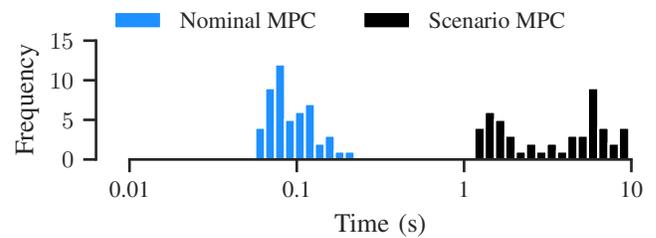}
    \caption{Histograms of time taken to solve the first MPC optimization for both nominal MPC
and scenario MPC.}
    \label{fig:results_time_histogram}
\end{figure}

Figure \ref{fig:results_time_histogram} shows the time taken for the first MPC optimization (i.e. the longest horizon) of each journey using both nominal and scenario MPC. The averages for nominal and scenario MPC were 0.098 s and 4.8 s, so the scenario approach does not meet the real-time requirement for the hardware used for these experiments. However, the ADMM iteration can be accelerated using parallel processing, and it has been demonstrated that the computational time of a similar ADMM algorithm can be significantly reduced using a GPU \cite{Qureshi2019} (the algorithm was implemented in vectorized sequential Matlab code here; the Matlab \emph{gpuArray} data type was investigated but not found to provide a performance benefit). Therefore, it is possible that the algorithm could be accelerated to meet the real-time requirement when implemented on suitable hardware. Additionally, the overall state-of-charge trajectories in Figure~\ref{fig:results_distance_1} closely match those in Figure~\ref{fig:results_distance_2}, so it may be possible to accelerate the algorithm and reduce the occurrence of constraint violations without negatively affecting performance (in terms of increased fuel consumption) by using fewer scenarios for the long-horizon predictions, and more for the one-step-ahead predictions. A systematic investigation into these possibilities is left to future work, but nevertheless, the ADMM algorithm still provides a clear benefit for this problem relative to general purpose convex optimization software as the solution times with CVX were found to be far too slow for the scenario MPC to be implemented in closed loop.

\section{Conclusion}

This paper presents a data-based approach to model predictive control for PHEV energy management, where the predictions of future driver behaviour are made directly from recorded examples of drivers previously completing a route. Results from scenario optimization are used to determine a bound on the confidence that the one-step-ahead optimization will be feasible with a given probability, and an efficient ADMM algorithm is proposed for the solution of the associated scenario MPC optimization problem. The performance of the scenario MPC algorithm is investigated through simulations of a passenger vehicle completing the same route multiple times, where it is found to obtain an extremely close approximation of the control inputs obtained using an MPC controller with full preview of future driver behaviour. 

The presented results have significant implications for future development and investigation of power-split optimal control techniques for PHEVs. Firstly, it has been demonstrated that previously recorded examples of a given route being driven are a powerful surrogate for the model-based prediction methods commonly found in the literature, and that long-horizon predictions of future driver behaviour can reliably be made directly from this data. Future work will investigate how well the performance observed in this study translates to more diverse driving conditions. Secondly, the results strongly suggest that future investigations on this topic should consider formulating the energy management problem in terms of \textit{completed distance}, instead of the time-based formulation commonly proposed in the literature.

\section*{ACKNOWLEDGMENT}

We would like to thank Professor Luigi del Re and
Dr Philipp Polterauer at Johannes Kepler University Linz, Austria,
for generously sharing the driver data used in this publication
and assisting with its processing.

\bibliographystyle{IEEEtran}
\bibliography{bibliography}

\begin{IEEEbiographynophoto}{Sebastian East}
received the M.Eng. degree in mechanical engineering from the University of Bath, Bath, U.K., in 2015, and the D.Phil. degree in engineering science from the University of Oxford, Oxford, U.K. in 2021. He is currently a Lecturer at the University of Bristol, U.K., where his research is focussed on problems in control, optimization, and machine learning.
\end{IEEEbiographynophoto}

\begin{IEEEbiographynophoto}{Mark Cannon}
received the M.Eng. degree in engineering science and the D.Phil. degree in control engineering from the University of Oxford, Oxford, U.K., in 1993 and 1998, respectively, and the M.S. degree in mechanical engineering from the Massachusetts Institute of Technology, Cambridge,
MA, USA, in 1995. He is currently an Associate Professor of engineering science with the University of Oxford and an Fellow of St John’s College, Oxford. His current research interests include robust and optimal control for constrained and uncertain systems, optimization for receding horizon control with robust constraints and stochastic uncertainty, and stochastic model
predictive control.
\end{IEEEbiographynophoto}

\end{document}